\documentclass[12pt]{report}
\usepackage{graphicx}
\usepackage{tikz}
\usepackage[utf8]{inputenc}
\usepackage[english]{babel} 
\usepackage{amsthm}
\usepackage{amsfonts}
\usepackage{array}
\usepackage{multirow}
\usepackage{enumitem}
\usepackage[fleqn]{amsmath}
\usepackage{amssymb}
\usepackage{thmtools}
\usepackage{dsfont}
\usepackage{rotating}

\theoremstyle{definition}

\declaretheoremstyle[
spaceabove=6pt, spacebelow=6pt,
headfont=\normalfont\bfseries,
notefont=\mdseries, notebraces={(}{)},
bodyfont=\normalfont,
postheadspace=1em
]{pbstyle}
\declaretheoremstyle[
spaceabove=6pt, spacebelow=13pt,
headfont=\normalfont\bfseries,
notefont=\mdseries, notebraces={(}{)},
bodyfont=\normalfont,
postheadspace=1em,
headpunct={:},
qed=$\blacktriangleleft$,
numbered=no
]{solstyle}

\usetikzlibrary{arrows}

\begin{document}
\begin{center}

{\LARGE \bf  The edge-vertex  inequality in   a planar graph  and a bipartition  for  the class of all  planar graphs}\\

\vspace*{.5in}

  {\bf  M. R.  Emamy-K.},\footnote{Dept. of Math.  Nat. Sc. UPR Rio Piedras, PR} \, {\bf Bahman  Kalantari},\footnote{Dept. of CS.,  Rutgers Univ. New Brunswick, NJ} \, {\bf Tatiana Correa}
\footnote {Nat. Sc. UPR Rio Piedras, PR} \\

\vspace*{.25in}

\end{center}

\vspace*{1in}
\begin{center}

{\LARGE  \bf Abstract}\\
\end{center}

\vspace*{.4in}

  For a   planar graph with a given f-vector   $(f_{0}, f_{1}, f_{2}),$ we introduce
 a cubic polynomial  whose coefficients depend on   the  f-vector. The planar graph is said to be real if all  the roots of  the corresponding   polynomial   are  real.  Thus we  have  a bipartition of all planar graphs
 into  two disjoint class of graphs, real and complex ones.
 As a contribution toward a full recognition of  planar graphs in this bipartition, we study   and recognize completely  a  subclass of
 planar graphs that includes all   the connected grid subgraphs. Finally, all the   2-connected triangle-free complex planar graphs of 7 vertices  are listed.

\newpage

  The inequality $ 3 (f_{0} - 2) \geq  f_{1} $ is an old  basic  relation  that holds  for all planar graphs, where its number of vertices $ f_{0} \geq 3$ and   $f_{1} $   stands for its number of edges, see  [1,2]. In other words, for   a fixed number of vertices $f_{0}$ a planar graph can not have more than  $ 3 (f_{0} - 2)$   edges.  Here is  a natural question: Is there  a  tighter  upper bound   for  $f_{1} $  that holds for  all the planar graphs?  The answer is; not very likely, due to the fact that  maximal planar triangulations have exactly  $ 3 (f_{0} - 2)$ number of edges. However,  this inequality
  has   a stronger  version; that is,   $ 2(f_{0} - 2) \geq  f_{1}, $  and    holds   for  triangle-free  planar graphs. The  fact that  the second inequality
  is  sharper  and  holds over a   subclass of planar graphs is leading  us to the idea behind this note.  Here, we introduce a  different inequality that is also stronger than the first one
 over  a subclass of planar graphs. This new inequality is  linked to  zeros of certain  Euler type  polynomials, and  leads us  to search for  a new bipartition   of  planar graphs.\\

Let  $G$  be a   simple  planar graph with no loops or parallel edges with   the usual  f-vector  $(f_{0}, f_{1}, f_{2}), $  where  $f_{0}, f_{1}, f_{2}$  stand  for the number of vertices, edges  and faces of $G.$ Consider the cubic polynomial $$p(x) =  f_{2}x^{3}+ f_{1}x^{2}+ f_{0}x + 2  $$ associated with  $G.$ Here   in this note,  $p(x)$  will be called Euler polynomial of $G.$  Evidently,  the   Euler polynomial  is uniquely defined by $G,$  however many non-isomorphic planar  graphs
may  have the same f-vector  and  hence the same Euler polynomial.   It is  rather trivial  to show that   the roots of $p(x)$ are all real if and only  if the relation\\
$$  (f_{0} + 2 )^{2}  \geq  8(f_{1} + 2 )\, \, \, \,  \, \,  (1) $$ \\

holds.
Our final    aim  here     is to  classify those   planar graphs that satisfy this  real root property. In other words, we  like  to  construct a  bipartition
 of all planar graphs into two real   and complex  subclasses.  The connection between real root properties of this type of polynomials   and  the associated  class of graphs  could open a link between graphs and  polynomial root finding. For more on root-finding and its application in polynomiography, another  facet  of mathematical arts,  see [3].
A planar graph that satisfy  the  inequality above is defined to be    a real graph,   and is called a complex  graph  otherwise.
 The inequality  can also be stated as:\\    $$  (f_{0} - 2 )^{2}  \geq   8f_{2} \, \, \, \,  \, \,   (2),$$  applying  the Euler relation for planar graphs.\\

It is intuitively  evident   that  the planar graphs with large number of vertices are real and on the other side the small graphs are complex.
The specific question is what do these large and small mean? and  what about  those graphs  in the middle levels? The main results of this note, Theorems 3-5, are toward  a   complete characterization of those real planar graphs of the middle levels.
It is easy to check that  all  the  planar graphs with 4 vertices  or less   are complex. The only planar graphs with 5 vertices that
are real must have $f_2 = 1,$  so  they are  acyclic that includes all the trees. Similarly,  the only real graphs of  6 vertices
are  the one  with at most one cycle.   It is  also easy to show that all the planar graphs with $  f_{0}  \geq 18$ are real. Hence, the middle level means those planar graphs with $ 7  \leq  f_{0}  \leq 17.$ If the  working  domain   is  restricted to  the triangle free ones,  then  the middle  level shrinks down to  $ 7  \leq  f_{0}  \leq 9.$       From now  on,  we assume  all the  graphs  are connected  unless otherwise is  stated.
 In the following  theorem, we begin with some small  class of  planar graphs where  this real property recognition is straight forward.

\noindent Let $P_n$  denote a path of size  $n$, that is  $P_n$   has $n $
vertices.   A rectangular  grid graph  (or simply a grid) $G_{m,n}$,  $m  \leq n$ is   defined to be the cartesian product
$P_{m}  \,  \times \, P_{n},$   that has $mn$  vertices,  $2mn -m-n$  edges   and $(m - 1 ) (n - 1 ) + 1 $ faces.\\

\noindent  { \bf Theorem  1.}

\noindent   (a)  A tree is real iff $  f_{0}  \geq  5$\\

 \noindent  (b) A cycle  is real iff $  f_{0}  \geq  6$\\

\noindent   (c)    The square grid  $G_n = P_{n}  \,  \times  \, P_{n}$ is real iff  $n \geq 3.$\\

\newpage

\noindent { \bf Proof.}  The statements of (a)  and  (b)  can directly  be verified. To prove (c),  applying the values of $f_{0}$ and  $f_{2} $  for the square  grid  of $n^2$ vertices,  condition (2) will be  stated as

$$  (n^2 - 2 )^{2}  \geq  8( n^2 -2n +2) \, \, \, \,  \, \, (3)   $$\\   Let us  define $f(n) = (n^2 - 2 )^{2}  -  8( n^2 -2n +2). $ Verifying  relations  $f(2)  < 0, $  $ f(3)>0 $ and   the fact that $f $   is  increasing for $n \geq 3$  shows that  the  inequality (3) holds,   and  this   completes the  proof.\\

  \noindent  { \bf Theorem 2.} A rectangular  grid  $G_{mn}$,    $  m,n  \geq  1$  is  always real except the following cases:\\

 (a)  $G_{1,n}$  for  $1 \leq  n  \leq 4$\\

  (b)  $G_{2,n}$  for  $2\leq  n  \leq 3$\\

 \noindent {\bf  Proof.}  A direct application of the definition shows  that $G_{m,n}$ is  complex
for all natural numbers $m,n, $ whenever its  number of vertices
  $mn\leq  4.  $   This   includes all the  paths with 4  vertices  or  less, and   a single   square.
   The grid  with   $mn = 5,  $  i.e. $G_{1,5}$, is  obviously real. The  case of    $mn = 6 $ includes  $G_{1,6}$ that is real  and $G_{2,3},$ the union of two adjacent  squares, that is complex.
 To show  that all  the cases  with   $mn \geq  7  $   are real,
   consider the inequality  that is deduced from (2), that is;\\

$$  (mn - 2 )^{2}  \geq  8( mn -m  -n +2) \, \, \, \, $$  or \\

$$  (mn - 2 )^{2} -  8mn \geq  16 - 8(m + n ) \, \, \, \, \, \, \, \,  (4)$$

\newpage

To complete  the proof,  note that the function $  g(x) =  (x - 2 )^{2} -  8x$
 is increasing in $x$  for $x \geq 7.$  And then    the relation (4) holds for all $mn \geq  7, $  since
 $16 - 8(m + n) < g(7) =  -31$  for all $mn \geq  7. $ \\

 The   real-complex  recognition problem  for  more    general planar graphs  will be  the subject of our  work
 in the  following theorems.
We  first  restate   the original basic  inequalities as the next lemma.\\

 \noindent  { \bf Lemma 1.}  Let G be any planar graph with $f_{0}$ and  $f_{1} $ vertices and  edges respectively,
 where $f_{0} \geq 3. $  Then,\\

\noindent  (a)   $ 3 (f_{0} - 2) \geq  f_{1} $\\

 \noindent (b)   Moreover   if $G$  is triangle- free, then $ 2(f_{0} - 2) \geq  f_{1} $\\

 We proceed by raising  the  critical question on the initial motivation  behind the note.
  How interesting are these  real or complex classes of planar graphs? This  depends on the size of the real  classes for $  5   \leq    f_{0}  \leq  17.  $   A large class of real graphs means that there is a stronger  inequality    than the linear one of Lemma 1-a  over  that  special   real class.

  The inequality   (a) of the lemma presents    $ 3 (f_{0} - 2)  $   as  a linear upper  bound for  $  f_{1}, $  however
 the upper  bound  for $  f_{1}$  in (1),  $ \frac { (f_{0} + 2 )^{2}}{8} - 2$,    is quadratic and in general   will  not  be superior  to the  linear one   for graphs with very  large number of vertices.
 In fact the quadratic inequality  would be  interesting exactly  when\\  $$ \frac { (f_{0} + 2 )^{2}}{8} - 2  \leq   3 (f_{0} - 2),$$\\  comparing the upper bounds   for  $  f_{1}$ in the two inequalities.  The latter implies $   f_{0}  \leq  17.   $  A similar comparison   for  the triangle free planar graphs
 shows  $   f_{0}  \leq  9.   $  The  existence of  these   bounds  form the first half of   the next theorem.\\
 
 \newpage

 \noindent  { \bf Theorem 3.}

 \noindent   Any planar graph with $f_{0} \geq  18$ is real, in addition, any triangle-free  planar graph with $f_{0} \geq  10$ is also real.
 That is:
 $$  (f_{0} - 2 )^{2}  \geq  8f_{2} \, \, \, \,      $$  holds for  planar graphs with $f_{0} \geq  18$, and for  any triangle-free  planar graph with $f_{0} \geq  10$. In  both cases these  are  the best possible inequalities,  in other words there are appropriate complex planar graphs  of 17 and 9 vertices  respectively.\\

 \noindent { \bf Proof.} The discussion above  concludes   that  $  (f_{0} + 2 )^{2}  \geq  8(3 f_{0} - 4 )$   for all  $f_{0}   \geq  18.$
 Then by applying Lemma 1 (a),   $3 f_{0} - 4  \geq f_{1} + 2,$ we  obtain:
$$  (f_{0} + 2 )^{2}  \geq  8(f_{1} + 2 )\, \, \, \,       $$  For the second part, again  it is  easy to see that
 $  (f_{0} + 2 )^{2}  \geq  16( f_{0} - 1 )$  for all  $f_{0}   \geq  10.$
 Then by applying Lemma 1 (b),   $2 ( f_{0} - 1)  \geq  f_{1} + 2,$ we  obtain:
 $$  (f_{0} + 2 )^{2}  \geq  8(f_{1} + 2 )\, \, \, \,       $$  or equivalently
 $$  (f_{0} - 2 )^{2}  \geq  8f_{2} \, \, \, \,  \, \, \, \, \, \, \, \,\, \, \, \,     $$

 To complete the proof, maximal  complex planar graphs  of 17 and  maximal  complex triangle-free planar graphs  of 9 vertices   are presented in Figures 1  and  2  respectively.\\

  The complex planar graph of
 17 vertices   above   is  a maximal  triangulation,  that is to say every  face including the unbounded face is a 3-cycle.  In fact  there are more
  maximal  triangulations that are  complex.

 { \bf Corollary.}  Any maximal planar triangulation graph  with   $  3   \leq    f_{0}  \leq  17  $  is complex, and any  complex
   planar graphs of 17  vertices is either a maximal triangulation   or  can be  obtained from such a  maximal  triangulation by removal of only  a  single  edge.
   
\newpage

 { \bf   Proof. }  Any maximal triangulation of $  f_{0} $  vertices  must have   $  f_{1} =3 f_{0} - 6 $  edges  and   $  f_{2} =2f_{0} - 4 $   faces.  So,
   to have a   complex maximal triangulation we  must have $(f_{0} -2)^ {2} < 8( 2f_{0} - 4).$ The latter is equivalent to $  3   \leq    f_{0}  \leq  17.  $
   Since  any  maximal triangulation of $  f_{0} = 17$  vertices has 45 edges (with 30 faces), and
    any complex planar graph of 17 vertices must have at least 44 edges (and  29 faces), then  the  last part of the corollary is concluded.



\begin{center}
\begin{tikzpicture}
\draw [gray!90!]  (-7,0)--(0,7)--(7,0) --(-7,0);
\draw   [gray!90!] (-3.5,1)--(3.5,1);
\draw  [gray!90!]  (-7,0)--(-3.5,1);
\draw  [gray!90!] (-7,0)--(-1.5,1);
\draw   [gray!90!] (-1,0)--(-1.5,1);
\draw   [gray!90!](-1,0)--(1,1);
\draw    [gray!90!](1,1)--(7,0);
\draw   [gray!90!] (7,0)--(3.5,1);
\draw    [gray!90!] (-1.5,1)--(-1,1.8);
\draw   [gray!90!](1,1)--(-1.1,1.8);
\draw   [gray!90!](-3,2.6)--(3,2.6);
\draw   [gray!90!](-3,2.6)--(-3.5,1);
\draw   [gray!90!](-3,2.6)--(-1,1.8);
\draw   [gray!90!](-1,1.8)--(-1.5,2.6);
\draw   [gray!90!](-1,1.8)--(1,2.6);
\draw   [gray!90!](1,2.6)--(3.5,1);
\draw   [gray!90!](3,2.6)--(3.5,1);
\draw   [gray!90!](1,2.6)--(1,1);
\draw   [gray!90!](7,0)--(3,2.6);
\draw   [gray!90!](-7,0)--(-3,2.6);
\draw   [gray!90!](-1.5,1)--(-3,2.6);
\draw   [gray!90!](-1.5,2.6)--(-.8,3.8);
\draw   [gray!90!](-.8,3.8)--(1,2.6);
\draw   [gray!90!](-.8,3.8)--(.8,3.8);
\draw   [gray!90!] (1,2.6)--(.8,3.8);
\draw   [gray!90!](.8,3.8)--(3,2.6);
\draw   [gray!90!](-3,2.6)-- (-.9,5);
\draw   [gray!90!](-3,2.6)-- (0,7);
\draw    [gray!90!](-.8,3.8)--  (-.9,5);
\draw    [gray!90!](-1.5,2.6)--  (-.9,5);
\draw    [gray!90!] (-.9,5)--  (0,7);
\draw    [gray!90!] (.8,3.8) --  (0,7);
\draw     [gray!90!]   (-.8,3.8)  --  (0,7);
\draw     [gray!90!]   (.9,5) --  (0,7);
\draw     [gray!90!]   (.9,5) -- (3,2.6);
\draw     [gray!90!]   (.9,5) -- (.8,3.8);
\draw       [gray!90!](3,2.6) --  (0,7);
\filldraw      (-7,0) node[below =2pt]{\bf{15}} circle (2.5pt)
                   (-1,0) node[below=2pt]{\bf{16}} circle (2.5pt)
                   (0,7) node[above=2pt]{\bf{1}} circle (2.5pt)
                   (7,0) node[below=2pt]{\bf{17}}circle (2.5pt)
                   (-3.5,1) node[right =2pt]{\bf{11}} circle (2.5pt)
                   (-1.5,1) node[right=3pt]{\bf{12}} circle (2.5pt)
                   (1,1) node[right=3pt]{\bf{13}} circle (2.5pt)
                   (3.5,1) node[right=3pt]{\bf{14}}circle (2.5pt)
		  (-1,1.8) node[right=3pt]{\bf{10}}circle (2.5pt)
		  (-3,2.6) node[above right=.1pt]{\bf{6}}circle (2.5pt)
		  (-1.5,2.6) node[above right=.1pt]{\bf{7}}circle (2.5pt)
		  (1,2.6) node[above right=.1pt]{\bf{8}}circle (2.5pt)
		  (3,2.6) node[above right=.1pt]{\bf{9}}circle (2.5pt)
		   (-.8,3.8) node[above right=.1pt]{\bf{4}}circle (2.5pt)
		    (.8,3.8) node[above right=.1pt]{\bf{5}}circle (2.5pt)
		    (-.9,5) node[above right=.1pt]{\bf{2}}circle (2.5pt)
		    (.9,5) node[above right=.1pt]{\bf{3}}circle (2.5pt);

 \draw (0, -1.5) node {Fig. 1};

 \end{tikzpicture}
 \end{center}

\newpage

\begin{center}
\begin{tikzpicture}[scale=.8]
      \filldraw [black] (0,0) circle (2pt);
       \draw (0,-.3) node{$a$};
      \filldraw [black] (0,5) circle (2pt);
   	 \draw (.4,5.1) node{$h$};

     \foreach \x in {-7,-3,3,7}
       {
     	\draw [gray] (\x,2) -- (0,5);
         }
	\foreach \x in {-7, -3,3 ,7}
	{
   \filldraw [black] (-5,2) circle (2pt);
   \filldraw [black] (5,2) circle (2pt);
    \filldraw [black] (\x,2) circle (2pt);
	  \draw [gray] (\x,2) -- (0,0);
     	}
     	\draw (-7.4,2) node{$b$};
     	\draw (7.2,2) node{$g$};
     	\draw (-4.7,2.2) node{$c$};
     	\draw (5.2,2.2) node{$f$};
     	\draw (-2.4,2) node{$d$};
     	\draw (3.2,2.2) node{$e$};
	\draw[gray!90!] (-7,2)--(-3,2);
	\draw [gray!90!] (7,2)--(3,2);
	\filldraw[black] (0,7) circle (2pt);
	\draw (.3,7.3) node{$k$};
	\draw[gray!90!] (0,7)--(-7,2);
	\draw[gray!90!] (0,7)--(7,2);

 \draw (0, -1.5) node {Fig. 2};

\end{tikzpicture}
\end{center}

\vspace*{.5in}

The complex triangle-free planar graph of  9 vertices  above  can be used to construct more of such examples with 8 and then again with 7 vertices.
This can be done by removal of a degree 2 vertex, applying the following lemma. \\

  {\bf    Lemma 2. }   Let   $G$  be a complex  planar  graph with $f_{0} \geq  7, $  and  a vertex $x$  of degree 2.    Then $G - x $  is also a complex  planar graph.

{ \bf  Proof.   }  It is   a direct consequence   of  the definition. \\

\newpage

 The last corollary shows that  only for $f_0  =  17$ there is  a   large class of real graphs and we show more of this in the next theorem.
    By a  grid  graph, we mean any
  graph that  can be embedded in  the square grid $G_{n}$ as a subgraph  for  a  large enough $n.$
 The  full characterization of real or complex planar graphs for all  $  7   \leq    f_{0}  \leq  17,  $  remains  to be an interesting computational problem. However,   we proceed to  show a  complete characterization of real  or complex subgraphs  of a grid  in the next theorem.\\


\noindent { \bf Theorem 4.} Any  connected  grid  graph with the number of vertices  $f_0  \geq  7$ is real.\\

\noindent { \bf Proof.}  First, consider only  the case  of 7 vertices,  and  suppose there is a  connected   complex grid  graph  $G$  with  $f_0  =  7.$  Evidently $G$   is planar, bipartite with at least 4 faces,
applying  $  (f_{0} - 2 )^{2} <  8f_{2}.$
 Three  of these  faces,  say   $C_{1},  C_{2},  C_{3},$  are cycles.
Since $G $ is bipartite   and     $f_{0} = 7,$   then,   these cycles  either must      be  4-cycles ( unit squares)  or 6-cycles  (rectangular 6-gons).  However  none of them can be a 6-cycle,
otherwise there would not be  any possibility   for  a  second cycle, by the restriction of $f_{0} = 7$ . Now, let   them  all to  be  4-cycles.  Since  the  union of any  two cycles has    at least  6 vertices
then  $  C_{1}\bigcup C_{2}\bigcup C_{3}$   needs   at least 8 vertices, a contradiction by $f_{0} = 7$ . The latter   will be apparent  if we   notice that the  connected  grid graph
$  C_{1}\bigcup C_{2}\bigcup C_{3}$   has only 5 options of non-isomorphic graphs, given that all the  cycles are   unit squares.  On the other hand,  if $  C_{1}\bigcup C_{2}\bigcup C_{3}$ is disconnected,  then  it requires even more  than 8 vertices,  that is
impossible. Finally, for $f_0  =  8$  and $f_0  =  9$ the minimum number of faces are 5  and  7 respectively. Then the  minimum number of bounded    faces  (  i.e., 4-cycles)
that are needed will be 4  and 6 respectively. We already needed  at least  8 vertices for only 3 cycles, so  the same contradiction is evident here.\\

The  grid  graphs are triangle-free,   however the  last theorem does not hold for arbitrary triangle-free planar graphs.  In fact,
 we have    presented an  example  of  a  complex  triangle-free planar graph of  9   vertices  in Figure 2, and  also examples of such complex planar graphs with 8 and 7 vertices, applying Lemma 2.  For  $f_{0}   =   7$ all examples of complex  triangle-free planar graphs,
   the   only possible number of edges are 9  and  10, applying the inequality (1) and Lemma 1 (b).  Hence,  the  number of edges
  of such  minimal  examples is  9,  and  this  is the motivation for  the next classification  theorem  that covers  the  corresponding  class of those examples . \\

 In the following theorem, let $P$ be a polygon and $S$ a finite set of points in the interior of $P$.  We  will explore all the possibilities for construction of  a  planar graph $G$ whose vertex set  $(vert P)\,  \cup \, S $ is fixed,  its f-vector is $(f_{0},f_{1},f_{2})=(7,9,4)$ and  includes the polygon $P$ as a cycle. In this construction,  those   edges that connect vertices of $S$ to  the vertices  of $P$  are  called connecting edges.
 In fact, we aim to  characterize all the non-isomorphic 2-connected triangle-free planar graphs  that    satisfy:
 $f=(f_{0},f_{1},f_{2})=(7,9,4)$, and  evidently all are   complex   graphs.
Since all the faces of a 2-connected planar graph are cycles, we may consider $G$ to  be a  2-connected triangle free planar graph with  $(f_{0},f_{1},f_{2})=(7,9,4)$ whose unbounded face is denoted by the  cycle $P$.\\

\noindent { \bf Theorem 5.}  All the  minimal  complex  non-isomorphic  2- connected triangle-free planar graphs with 7 vertices are listed in the  figures  3, 4  and 5.

\noindent {\bf Proof}
 We may assume that such graph $G$ has no leaf, or any other cut vertices, since it is 2-connected. Moreover, all of its faces are cycles.
Consider a planar drawing of $G$, and let $P$  be the outer face  of this drawing, that is a convex polygon that includes all other vertices  in its interior. Obviously,  $P$ can be a   hexagon with   only  one  vertex left inside  the hexagon, a pentagon with two vertices inside or a quadrilateral with three vertices inside. Of course,   the condition  $f_{2}=4$,$f_{1}=9$ and $f_{0}=7$  holds  in all the 3 cases.

Consequently,  the characterization is divided into three cases:
 	Case 1: $P$ is a 6-cycle and $\mid S \mid=1$.
	Case 2: $P$ is a 5 -cycle with $\mid S \mid =2$
	 and in Case 3: $P$ is a  4- cycle with $\mid S \mid =3$

Figure 3 shows the unique graph $G$, where the three connecting edges of this figure are the only possibility since triangular faces are not allowed. The Case 2 includes 4 graphs where the two internal vertices in $S$ have degree 2 or 3, by the 2- connectivity of $G$. In the first two graphs  Figure 4-a, there is one  edge incident to the two vertices of $S$, and so there are $3$ connecting edges. In the graphs of Figure 4-b there is no internal edge and so there are $4$ connecting edges.

\newpage

	In Figure 4-a, again triangles are not allowed, so the two connecting edges that are incident to a vertex of $S$ can only form a 4-cycle or a 5-cycle joint with more edges from the cycle $P$. In Figure 4-b, two connecting edges incident to a vertex in $S$  must also  be part of a  4-cycle, or a 5-cycle joint with more edges from the cycle $P,$  since triangles are forbidden.  The two remaining  connecting edges are incident to the other vertex of $S.$  Hence, there will remain   only two possibilities  left in Figure  4-b.

 Case 3 includes 7 more graphs, where in all of them $|S|=3$. However in the first five graphs, Figure 5-a,  there are two edges incident to vertices of $S$, and 3 more connecting edges. In the last two graphs, Figure 5-b,  there is only one edge incident to the  vertices of $S$  and  four connecting edges.

\begin{tabular}{c c}

 \begin{tikzpicture}[scale=.8]
    \draw [gray!90!]  (0,0)--(4,0);
    \draw [gray!90!]  (0,0)--(-2,2);
    \draw [gray!90!]  (-2,2)--(0,4);
    \draw [gray!90!]  (0,4)--(4,4);
    \draw [gray!90!]  (4,4)--(6,2);
    \draw [gray!90!]  (6,2)--(4,0);
    \draw [gray!90!]  (4,0)--(2,2);
        \draw [gray!90!]  (2,2)--(4,4);
         \draw [gray!90!]  (-2,2)--(2,2);
\filldraw
                   (0,0) node[below =2pt]{\bf{2}} circle (2.5pt)
                   (-2,2) node[below =2pt]{\bf{1}} circle (2.5pt)
                   (4,4) node[below =2pt]{\bf{1}} circle (2.5pt)
                   (6,2) node[below =2pt]{\bf{2}} circle (2.5pt)
                   (4,0) node[below =2pt]{\bf{1}} circle (2.5pt)
 (0,4) node[below =2pt]{\bf{2}} circle (2.5pt)
 (2,2) node[below =2pt]{\bf{2}} circle (2.5pt);

 \draw (2,-2) node {Fig. 3};


\end{tikzpicture}

\\
\\
 \\\\\\
%


    \begin{tikzpicture}[scale=.8]

    \draw [gray!90!]  (-3,1)--(0,3);
    \draw [gray!90!]  (-3,1)--(-2,-2);
    \draw [gray!90!]  (-2,-2)--(-1.5,1);
    \draw [gray!90!]  (-1.5,1)--(0,3);
    \draw [gray!90!]  (-1.5,1)--(1.5,-.75);
    \draw [gray!90!]  (1.5,-.75)--(3,1);
    \draw [gray!90!]  (2,-2)--(3,1);
        \draw [gray!90!]  (0,3)--(3,1);
        \draw [gray!90!] (-2,-2)--(2,-2);
\filldraw      (-3,1) node[above =2pt]{\bf{}} circle (2.5pt)
                   (-2,-2) node[below =2pt]{\bf{}} circle (2.5pt)
                   (2,-2) node[below =2pt]{\bf{}} circle (2.5pt)
                   (-1.5,1) node[below =2pt]{\bf{}} circle (2.5pt)
                   (0,3) node[above =2pt]{\bf{}} circle (2.5pt)
                   (1.5,-.75) node[below =2pt]{\bf{}} circle (2.5pt)
                  ( 3,1) node[above =2pt]{\bf{}} circle (2.5pt);

\end{tikzpicture}


 \hspace{30pt}

    \begin{tikzpicture}[scale=.8]

    \draw [gray!90!]  (-3,1)--(0,3);
    \draw [gray!90!]  (-3,1)--(-2,-2);
    \draw [gray!90!]  (-3,1)--(-1.5,0);
    \draw [gray!90!]  (-1.5,0)--(0.5,1);
    \draw [gray!90!]  (-2,-2)--(0.5,1);
    \draw [gray!90!]  (2,-2)--(3,1);
    \draw [gray!90!]  (0.5,1)--(0,3);
        \draw [gray!90!]  (0,3)--(3,1);
        \draw [gray!90!] (-2,-2)--(2,-2);
\filldraw      (-3,1) node[above =2pt]{\bf{}} circle (2.5pt)
                   (-2,-2) node[below =2pt]{\bf{}} circle (2.5pt)
                   (2,-2) node[below =2pt]{\bf{}} circle (2.5pt)
                   (-1.5,0) node[below =2pt]{\bf{}} circle (2.5pt)
                   (0,3) node[above =2pt]{\bf{}} circle (2.5pt)
                   (0.5,1) node[below =2pt]{\bf{}} circle (2.5pt)
                  ( 3,1) node[above =2pt]{\bf{}} circle (2.5pt);

     \end{tikzpicture}

\\
\\

     \begin{tikzpicture}[scale=.8]

      \draw (-10,0) node {Fig. 4-a};
      \end{tikzpicture}
\\
\\
\\
\\


   \begin{tikzpicture}[scale=.8]
    \draw [gray!90!]  (-3,1)--(0,3);
    \draw [gray!90!]  (-3,1)--(-2,-2);
    \draw [gray!90!]  (-2,-2)--(-1,1);
    \draw [gray!90!]  (-1,1)--(0,3);
    \draw [gray!90!]  (2,-2)--(3,1);
    \draw [gray!90!]  (2,-2)--(1,1);
    \draw [gray!90!]  (1,1)--(0,3);
        \draw [gray!90!]  (0,3)--(3,1);
        \draw [gray!90!] (-2,-2)--(2,-2);
\filldraw      (-3,1) node[above =2pt]{\bf{}} circle (2.5pt)
                   (-2,-2) node[below =2pt]{\bf{}} circle (2.5pt)
                   (2,-2) node[below =2pt]{\bf{}} circle (2.5pt)
                   (-1,1) node[below =2pt]{\bf{}} circle (2.5pt)
                   (0,3) node[above =2pt]{\bf{}} circle (2.5pt)
                   (1,1) node[below =2pt]{\bf{}} circle (2.5pt)
                  ( 3,1) node[above =2pt]{\bf{}} circle (2.5pt);
\end{tikzpicture}

 \hspace{30pt}


    \begin{tikzpicture}[scale=.8]

    \draw [gray!90!]  (-3,1)--(0,3);
    \draw [gray!90!]  (-3,1)--(-2,-2);
    \draw [gray!90!]  (-2,-2)--(1,0.25);
    \draw [gray!90!]  (1,0.25)--(0,3);
    \draw [gray!90!]  (2,-2)--(3,1);
    \draw [gray!90!]  (-2,-2)--(2,0.25);
    \draw [gray!90!]  (2,0.25)--(0,3);
        \draw [gray!90!]  (0,3)--(3,1);
        \draw [gray!90!] (-2,-2)--(2,-2);
\filldraw      (-3,1) node[above =2pt]{\bf{}} circle (2.5pt)
                   (-2,-2) node[below =2pt]{\bf{}} circle (2.5pt)
                   (2,-2) node[below =2pt]{\bf{}} circle (2.5pt)
                   (1,0.25) node[below =2pt]{\bf{}} circle (2.5pt)
                   (0,3) node[above =2pt]{\bf{}} circle (2.5pt)
                   (2,0.25) node[below =2pt]{\bf{}} circle (2.5pt)
                  ( 3,1) node[above =2pt]{\bf{}} circle (2.5pt);

     \end{tikzpicture}

     \begin{tikzpicture}[scale=.8]

      \end{tikzpicture}

\\

     \begin{tikzpicture}[scale=.8]

      \draw (-10,0) node {Fig. 4-b};
      \end{tikzpicture}

 \end{tabular}

\begin{tabular}{c c}
    \begin{tikzpicture}
    \draw [gray!90!]  (0,0)--(3,0);
    \draw [gray!90!]  (0,0)--(0,3);
    \draw [gray!90!]  (3,0)--(3,3);
    \draw [gray!90!]  (0,3)--(3,3);
    \draw [gray!90!]  (3,3)--(1.5,2);
    \draw [gray!90!]  (1.5,2)--(1,1);
    \draw [gray!90!]  (1.5,2)--(2,1);
        \draw [gray!90!]  (0,0)--(1,1);
            \draw [gray!90!]  (3,0)--(2,1);
\filldraw      (3,0) node[below =2pt]{\bf{}} circle (2.5pt)
                   (0,0) node[below =2pt]{\bf{}} circle (2.5pt)
                   (0,3) node[above =2pt]{\bf{}} circle (2.5pt)
                   (3,3) node[above =2pt]{\bf{}} circle (2.5pt)
                   (1.5,2) node[below =2pt]{\bf{}} circle (2.5pt)
                   (1,1) node[below =2pt]{\bf{}} circle (2.5pt)
                   (2,1) node[below =2pt]{\bf{}} circle (2.5pt);
\end{tikzpicture}


\hspace{30pt}
    \begin{tikzpicture}
    \draw [gray!90!]  (0,0)--(3,0);
    \draw [gray!90!]  (0,0)--(0,3);
    \draw [gray!90!]  (3,0)--(3,3);
    \draw [gray!90!]  (0,3)--(3,3);
    \draw [gray!90!]  (3,3)--(1.7,2.1);
    \draw [gray!90!]  (1.7,2.1)--(1,1.5);
    \draw [gray!90!]  (1,1.5)--(0,3);
    \draw [gray!90!]  (2,1.1)--(3,0);
\draw [gray!90!]  (2,1.1)--(1.7,2.1);
\filldraw      (3,0) node[below =2pt]{\bf{2}} circle (2.5pt)
                   (0,0) node[below =2pt]{\bf{1}} circle (2.5pt)
                   (0,3) node[above =2pt]{\bf{2}} circle (2.5pt)
                   (3,3) node[above =2pt]{\bf{1}} circle (2.5pt)
                   (1.7,2.1) node[below =2pt]{\bf{2}} circle (2.5pt)
                   (1,1.5) node[below =2pt]{\bf{1}} circle (2.5pt)
                   (2,1.1) node[below =2pt]{\bf{1}} circle (2.5pt);
\end{tikzpicture}

  \hspace{30pt}


    \begin{tikzpicture}
    \draw [gray!90!]  (0,0)--(3,0);
    \draw [gray!90!]  (0,0)--(0,3);
    \draw [gray!90!]  (3,0)--(3,3);
    \draw [gray!90!]  (0,3)--(3,3);
    \draw [gray!90!]  (3,3)--(1.5,1.2);
    \draw [gray!90!]  (1.5,1.2)--(2.5,1.7);
 \draw [gray!90!]  (1.5,1.2)--(0,0);
    \draw [gray!90!]  (2.5,1.7)--(2.3,1);
   \draw [gray!90!]  (3,0)--(2.3,1);

\filldraw      (3,0) node[below =2pt]{\bf{}} circle (2.5pt)
                   (0,0) node[below =2pt]{\bf{}} circle (2.5pt)
                   (0,3) node[above =2pt]{\bf{}} circle (2.5pt)
                   (3,3) node[above =2pt]{\bf{}} circle (2.5pt)
                   (1.5,1.2) node[below =2pt]{\bf{}} circle (2.5pt)
                   (2.3,1) node[below =2pt]{\bf{}} circle (2.5pt)
                   (2.5,1.7) node[below =2pt]{\bf{}} circle (2.5pt);
\end{tikzpicture}
\\
\\
 \\\\\\

    \begin{tikzpicture}
      \draw [gray!90!]  (1,1)--(3,0);
    \draw [gray!90!]  (0,0)--(3,0);
    \draw [gray!90!]  (0,0)--(0,3);
    \draw [gray!90!]  (3,0)--(3,3);
    \draw [gray!90!]  (0,3)--(3,3);
    \draw [gray!90!]  (0,3)--(1,2);
    \draw [gray!90!]  (1,2)--(2,1.5);
    \draw [gray!90!]  (2,1.5)--(1,1);
\draw[gray!90!] (1,1)--(0,3);

\filldraw
                   (3,0) node[below =2pt]{\bf{1}} circle (2.5pt)
                   (0,0) node[below =2pt]{\bf{2}} circle (2.5pt)
                   (0,3) node[above =2pt]{\bf{1}} circle (2.5pt)
                   (3,3) node[above =2pt]{\bf{2}} circle (2.5pt)
                   (1,2) node[below =2pt]{\bf{2}} circle (2.5pt)
                   (1,1) node[below =2pt]{\bf{2}} circle (2.5pt)
                   (2,1.5) node[below =2pt]{\bf{1}} circle (2.5pt);
                  (1,1) node [below=2pt{\bf{2}} circle (2.5pt);

\end{tikzpicture}

\hspace{30pt}


    \begin{tikzpicture}
    \draw [gray!90!]  (0,0)--(3,0);
    \draw [gray!90!]  (0,0)--(0,3);
    \draw [gray!90!]  (3,0)--(3,3);
    \draw [gray!90!]  (0,3)--(3,3);
    \draw [gray!90!]  (0,3)--(1,1);
     \draw [gray!90!]  (1,1)--(0.9,2.1);
    \draw [gray!90!]  (0.9,2.1)--(1.5,2.3);
    \draw [gray!90!]  (1,1)--(3,3);
\draw[gray!90!] (1.5,2.3)--(0,3);

\filldraw
                   (3,0) node[below =2pt]{\bf{}} circle (2.5pt)
                   (0,0) node[below =2pt]{\bf{}} circle (2.5pt)
                   (0,3) node[above =2pt]{\bf{}} circle (2.5pt)
                   (3,3) node[above =2pt]{\bf{}} circle (2.5pt)
                   (0.9,2.1) node[above =2pt]{\bf{}} circle (2.5pt)
                   (1,1) node[below =2pt]{\bf{}} circle (2.5pt)
                   (1.5,2.3) node[below =2pt]{\bf{}} circle (2.5pt);

\end{tikzpicture}

\\

\\

 \begin{tikzpicture}[scale=.8]

     \draw (-10,-6) node {Fig. 5-a};

 \end{tikzpicture}


 \\\\\\
\\\\
    \begin{tikzpicture}
    \draw [gray!90!]  (0,0)--(3,0);
    \draw [gray!90!]  (0,0)--(0,3);
    \draw [gray!90!]  (3,0)--(3,3);
    \draw [gray!90!]  (0,3)--(3,3);
    \draw [gray!90!]  (0,3)--(.8,2.1);
    \draw [gray!90!]  (.8,2.1)--(1.5,1.5);
    \draw [gray!90!]  (1.5,1.5)--(3,3);
\draw[gray!90!] (1,1)--(3,0);
\draw[gray!90!] (1,1)--(0,3);

\filldraw      (3,0) node[below =2pt]{\bf{2}} circle (2.5pt)
                   (0,0) node[below =2pt]{\bf{1}} circle (2.5pt)
                   (0,3) node[above =2pt]{\bf{2}} circle (2.5pt)
                   (3,3) node[above =2pt]{\bf{1}} circle (2.5pt)
                   (.8,2.1) node[below =2pt]{\bf{1}} circle (2.5pt)
                   (1,1) node[below =2pt]{\bf{1}} circle (2.5pt)
                   (1.5,1.5) node[below =2pt]{\bf{2}} circle (2.5pt);
                  (1,1) node [below=2pt{\bf{2}} circle (2.5pt);
\end{tikzpicture}
\hspace{30pt}


    \begin{tikzpicture}
    \draw [gray!90!]  (0,0.5)--(3,0.5);
    \draw [gray!90!]  (0,0.5)--(0,3.5);
    \draw [gray!90!]  (3,0.5)--(3,3.5);
     \draw[gray!90!] (0,3.5)--(3,3.5);
    \draw [gray!90!]  (0,3.5)--(1.5,2.5);
    \draw [gray!90!]  (0,3.5)--(1.1,1.5);
    \draw [gray!90!]  (1.1,1.5)--(3,0.5);
    \draw [gray!90!]  (1.5,2.5)--(2, 2);
\draw[gray!90!] (2, 2)--(3,0.5);

\filldraw     (2, 2)  node[below =2pt]{\bf{}} circle (2.5pt)
                   (3,0.5) node[below =2pt]{\bf{}} circle (2.5pt)
                   (0,0.5) node[below =2pt]{\bf{}} circle (2.5pt)
                   (0,3.5) node[above =2pt]{\bf{}} circle (2.5pt)
                   (3,3.5) node[above =2pt]{\bf{}} circle (2.5pt)
                   (1.5,2.5) node[above =2pt]{\bf{}} circle (2.5pt)
                   (1.1,1.5) node[below =2pt]{\bf{}} circle (2.5pt);

                  (2,1.7) node[below=2pt{\bf{2}} circle(2.5pt);
\end{tikzpicture}

\\

\\

\begin{tikzpicture}[scale=.8]

     \draw (-10,0) node {Fig. 5-b};
      \end{tikzpicture}
\\\\\\\\\\

  \end{tabular}

\noindent { \bf Final remark    and  more open problems }\\

It is  easy to extend the result of Theorem 5 to non-minimal complex planar graphs of 7 vertices,  since   those planar graphs  must
either  have 9   or  10 edges. Therefore, those  planar graphs of 7 vertices  with  10  edges can be obtained  from the  ones of 9 edges
by adding   one  more edges. However, the new extra edge can not create a triangular face, that is, it must be  added into a
hexagonal face only.  Thus, there  are only the bipartite graphs of   Figures  3  and 5,  the  ones with the  labeled vertices, that should be selected.
The   complete characterization of real or complex graphs with $  7   \leq    f_{0}  \leq  17$ remains  to be an interesting problem, and a practical computational project.
The notion of real and complex planar graphs that are defined  here are   directly connected to the quadratic
inequality that defines them. So,  the first natural question will be   if  there exist other class of such inequalities.
Are there stronger relation that can provide    tighter upper bound  for $f_0$ in terms of $f_1$?
We have noticed the real class of planar graphs are closely connected to the real root  property of Euler polynomials. How about
search for higher degree polynomials  such as higher degree Euler polynomials for polytopal graphs? The latter type of polynomials
can define a new  class of graphs   and    a new  link  between graphs  and   polynomial root-finding, see [3].\\

\vspace*{1in}

\noindent { \bf References}\\

\noindent [1] Bondy  J.A.  and  Murty  U.S.R.,  Graph Theory  with  Applications.   Springer, North Holland,  New  York (2008).

\noindent [2]  West, D.  B.,   Introduction to Graph Theory. Prentice Hall, NJ  (2001).

  \noindent [3] Kalantari   B.,  Polynomial   Root-Finding  and  Polynomiography.  World  Scientific, NJ  (2009).

\end{document}